\documentclass[twoside]{amsart}
\usepackage{amsfonts,amssymb,amscd,verbatim, enumerate}

\input xy
\xyoption{all}
\SelectTips{cm}{} 
\newcommand{\lra}{\longrightarrow}

\newcommand\Equal{\ar@{=}}

\newcommand{\cala}{{\mathcal A}}

\newcommand{\calo}{{\mathcal O}}

\newcommand{\cals}{{\mathcal S}}

\newcommand{\C}{\mathbb C}
\newcommand{\HH}{\mathbb H}
\newcommand{\R}{\mathbb R}
\newcommand{\Q}{\mathbb Q}
\newcommand{\Z}{\mathbb Z}

\renewcommand{\P}{\mathbb P}
\renewcommand{\L}{\mathbb L}

\newtheorem{definition}{Definition}[section]

\newtheorem{remark}[definition]{Remark}

\newtheorem{example}[definition]{Example}
\newtheorem{examples}[definition]{Examples}

\newtheorem{conjecture}[definition]{Conjecture}

\newtheorem{proposition}[definition]{Proposition}
\newtheorem{theorem}[definition]{Theorem}
\newtheorem{lemma}[definition]{Lemma}
\newtheorem{corollary}[definition]{Corollary}

\newcommand{\inj}{\hookrightarrow}
\newcommand{\surj}{\twoheadrightarrow}
\newcommand{\tensor}{\otimes}

\begin{document}

\title{An Index Theorem for Modules on a Hypersurface Singularity}

\author{Ragnar-Olaf Buchweitz}
\address{Dept.\ of Computer and Math\-ematical Sciences,
University  of Tor\-onto at Scarborough, 
1265 Military Trail, 
Toronto, ON M1C 1A4,
Canada}
\email{ragnar@utsc.utoronto.ca}

\author{Duco van Straten}
\address{
Fachbereich 17, 
AG Algebraische Geometrie, 
Johannes Gutenberg-UniversitŠt, 
D-55099 Mainz,
Germany}
\email{vanstrat@uni-mainz.de}
\date{\today}

\maketitle
\phantom{} \hfill {\em To the memory of V.~I.~Arnol'd}

\begin{abstract}
A topological interpretation of Hochster's Theta pairing of two mo\-dules on
a hypersurface ring is given in terms of linking numbers. This generalizes
results of M.~Hochster and proves a conjecture of J.~Steenbrink. As a
corollary we get that the Theta pairing vanishes for isolated
hypersurface singularities in an odd number of variables, as was conjectured
by H.~Dao. 
\end{abstract}

\section*{Introduction}

The interplay between topology and algebra is a central theme
in singularity theory. The formula expressing the Milnor number of
an isolated hypersurface singularity in terms of the length of its local 
algebra and the Eisenbud-Levine theorem expressing the topological degree as the
index of a quadratic form on the local algebra are cases in point.
In this paper we prove a theorem in the same spirit: we will 
show that a purely algebraic quantity, {\em Hochster's Theta pairing}
associated to two modules on the ring of an isolated hypersurface
singularity, can be expressed as a {\em linking number} of two associated
cycles. In this respect it is reminiscent of the classical interpretation 
of the intersection multiplicity in terms of linking, due to Lefschetz. 
The main ideas of this theorem were developed about $25$ years ago, 
during a visit of the first author to Leiden, where he lectured on the 
triangulated structure of the stable category of maximal Cohen-Macaulay 
modules, \cite{buchweitz}. On that occasion, J.~Steenbrink came up with 
a surprising conjecture relating Hochster's Theta pairing to the 
variation mapping in the cohomology of the Milnor fibre.\\

The past years have seen renewed interest in maximal Cohen-Macaulay modules
and matrix factorisations after the relation with $D$-branes and mirror
symmetry in Landau-Ginsburg models was suggested by M. Kontsevich and established by 
D. Orlov, \cite{orlov}.
The physics literature on the subject is rather enormous by now. The Kapustin-Li formula 
\cite{kapustinli} implies still another algebraic expression for the Theta pairing, however the 
main result of this paper still seems to be new.\\ 
The structure of the paper is as follows. In the first section we review
basic notation and definitions that we need. In the second section we 
formulate our main theorem. In the third section we use higher algebraic $K$-theory to 
define, for a module $M$ on a hypersurface ring, a class $\{M\}$ in $K^1$ and reformulate 
Hochster's Theta pairing in $K$-theoretic terms. In the fourth section we map everything 
to topological $K$-theory and give a different description for $\{M\}_{top}$. The
vanishing of $\theta(M,N)$ for an {\em odd number of variables} follows from the 
fact that in that case the Milnor fibre has only {\em even} cohomology. This
confirms a conjecture formulated in \cite{dao} that was  proven in the case of
quasi-homogeneous singularities and graded modules in \cite{moorecs}, and recently 
established as well in our context by Polishchuk and Vaintrob \cite{polishchukvaintrob} from the
(DG--)categorical point of view. 
In the fifth section we use the Chern-character to map topological $K$-theory to 
cohomology, which then leads to a proof of our main result. 
In a way our result represents a (rather hard won) triumph of topology over 
algebra that Arnol{'}d might have appreciated.

\section{Preliminaries}

{\bf Topology of singularities.}
We consider an isolated hypersurface singularity
\[ f \in  P:=\C \{x_0,x_1,\ldots,x_n\},\;\;f(0)=0\]
and we will choose a {\em good representative} \cite {agv},\cite{looijenga} 
\[f:X \lra D\] for the defining function in the usual way: first we take a sufficiently small ball 
$B(0,\epsilon)\subset\C^{n+1}$
such that the spheres $\partial B(0,\epsilon')$, for $0 <\epsilon' \le \epsilon$, 
are transverse to $f^{-1}(0)$, then put $X:=B(0,\epsilon) \cap f^{-1}(D)$,
where $D:=D_{\eta}$ is the $\eta$-disc in $\C$. Here $\eta$ is taken so small
that all fibres  $X_t:=f^{-1}(t), t\in D$ are transverse to the {\em Milnor sphere} 
$S:=\partial B(0,\epsilon)$. 
The fibre $X_0$ has an isolated singularity in $0$ and the restriction
$f^*:X^*:=X\setminus X_0 \lra D^*:=D\setminus \{0\}$ of $f$ is a locally 
trivial fibre bundle, the {\em Milnor fibration}.
The Milnor-fibres $X_t:=f^{*-1}(t)$ have the homotopy type of a wedge of
$n$-spheres, where the number $\mu(f)$ of these can be computed algebraically as
\[\mu(f)=\dim_{\C} (P/J_f)\,,\]  
the quotient of the power series ring modulo  the jacobian ideal of $f$,
\[ 
J_f:=\left(\frac{\partial f}{\partial x_0},\frac{\partial f}{\partial x_1},\ldots,
\frac{\partial f}{\partial x_n}\right) \subseteq P\,.
\]
The singular fibre $X_0$ is homeomorphic to the cone over the manifold 
$L:=S \cap X_0=\partial X_0$, the {\em link of the singularity}, that  is homotopy equivalent 
to the complex manifold $U:=X_0\setminus \{0\}$. 
By the triviality of the Milnor fibration near the boundary, one can identify 
$L$ also with the boundary $\partial X_t$, for any $t\in D^*$. From the 
homology sequence 
\[ 
0 \lra H_n(L) \lra H_n(F) \lra H_n(F,\partial F) \lra H_{n-1}(L) \lra 0
\]
of the pair $(X_t,\partial X_t)$, where we have written $F=X_t$ for the typical Milnor fibre, 
we see that $L$ has non-trivial homology groups only in degree $n-1$ and $n$ and 
that these are put in duality by the intersection product. 
For more details we refer to \cite{agv} and \cite{looijenga}.\\

{\bf Modules over Hypersurface rings.}
We will consider (finitely generated) modules $M$ over the hypersurface ring $R:=P/(f)$, a local 
ring of dimension $n$. Typical modules arise by considering {\em subvarieties} $Z$ lying inside 
the singular fibre $X_0$. Algebraically, we consider the ideal $I\subset P$ of $Z$, and as $Z$ is 
supposed to be contained in $X_0$, one has $f \in I$. The subvariety $Z$ now determines a 
module 
\[ \calo_Z:=P/I\,.\]

It is a basic fact, discovered by Eisenbud \cite{eisenbud}, that $R$-modules have a minimal 
resolution that is {\em eventually $2$-periodic}. In fact, a choice of generators for $M$ defines a 
surjection of a free $R$-module $F$ onto $M$, with the syzygy module $syz(M)$ as kernel:
\[ 0 \lra syz(M) \lra F \lra M \lra 0\,.\]
It follows from the depth-lemma that $depth_R(syz(M)) > depth_R(M)$ as long as 
$depth_R(M) < n$. Repeating the procedure with $syz(M)$, we see that after $n$ steps
we have an exact sequence of the following form
\[ 0 \lra M' \lra F_{n-1} \lra F_{n-2}\lra \cdots F_0 \lra M \lra 0\,, \]
where the  $F_i$ are free $R$-modules of finite rank and $depth_R(M')=n$. If $M'=0$, then 
$M$ has a free resolution of finite length, if $M' \neq 0$, then $M'$ is a 
{\em maximal Cohen-Macaulay module}, that is, $depth_R(M')=n$. 
So, ``up to free modules'', any $R$-module $M$ can be ``replaced'' by a maximal 
Cohen-Macaulay module. For a systematic account of such {\em Cohen-Macaulay 
approximations} we refer to \cite{buchweitz}, \cite{ausbuch}.
If $M$ is a maximal Cohen-Macauly $R$-module that is minimally generated by $p$ elements, 
its resolution as $P$-module has the form
\[ 0 \lra P^{p} \stackrel{A}{\lra} P^{p} \lra M \lra 0\,,\]
where $A$ is some $p \times p$-matrix, with%
\footnote{at least when $f$ is irreducible. In general, $\det A$ will divide $f^{q}$, with the exponent of 
an irreducible factor of $f$ in $\det A$ equal to the rank of $M$ on the corresponding component; 
see \cite{eisenbud}.}
 $\det(A)=f^q$. Expressing the
fact that multiplication by $f$ acts as $0$ on $M$ produces a matrix
$B \in Mat(p\times p, P)$, as in the diagram
\[
\xymatrix{
0 \ar[r] &P^{p} \ar[r]^-{\displaystyle A}\ar[d]_-{\displaystyle f\cdot\,}& 
P^{p} \ar[r] \ar[d]^-{\displaystyle f\cdot\,}\ar[dl]_{\displaystyle B} & 
M \ar[r] \ar[d]^{\displaystyle 0}& 0\\
0 \ar[r] &P^{p} \ar[r]^-{\displaystyle A}& 
P^{p} \ar[r] & 
M \ar[r] & 0\\
}
\]
such that
\[ A \cdot B = B \cdot A = f\cdot I\,,\]
where $I$ is the identity matrix. In other words, we find a {\em matrix factorisation} $(A,B)$ 
of $f$, determined uniquely, up to base change in the free modules $P^{p}$, by $M$. 
This matrix factorisation not only determines $M$, as $M = Coker(A)$, but also determines a 
resolution of $M$ as $R$-module 
$$ \cdots \lra R^p \stackrel{A}{\lra} R^p \stackrel{B}{\lra} R^p 
\stackrel{A}{\lra} R^p \lra M \lra 0$$
that  is plainly $2$-periodic. So the minimal resolution of an $R$-module looks in
general as follows:
\[\]
\[ \cdots \lra G \stackrel{A}{\lra} F\stackrel{B}{\lra} G\stackrel{A}{\lra}F \lra F_{n-1} \lra \cdots \lra F_0
\lra M \lra 0, \]
where $F= G$ are free $R$-modules.
As a consequence, for modules $M$ and $N$
over hypersurface rings, all homological algebra invariants
like $Tor_k^R(M,N)$ and $Ext_R^k(M,N)$ are eventually $2$-periodic.
It will be convenient also to consider the so-called {\em complete resolution}
$C^{\bullet}(M)$ of $M$, that is the bi-infinite $2$-periodic complex
\[ C^{\bullet}(M):\;\;\;\;\;\cdots \stackrel{B}{\lra} G \stackrel{A}{\lra} F\stackrel{B}{\lra} 
G\stackrel{A}{\lra}F \stackrel{B}{\lra} \cdots, \]
where by convention we always will put $F$ in {\em even} and $G$ in {\em odd} spots, although 
always $F=G$. For details we refer to \cite{eisenbud}, \cite{yoshino}, \cite{buchweitz}.\\

{\bf Hochster's Theta pairing.} 
The following quantity was considered  by Hochster in \cite{hochster}. 
In \cite{moorecs} it is called {\em Hochster's Theta invariant}.\\

{\bf Definition:}
The {\em Theta pairing} of modules $M$ and $N$ over the hypersurface ring 
$R=P/(f)$ is
\[
\theta(M,N):=length(Tor_{even}^R(M,N))-length(Tor_{odd}^R(M,N))\,,
\]
where 
\[Tor_{even}^R(M,N):=Tor_{2k}^R(M,N),\;\;\; Tor_{odd}^R(M,N):=Tor_{2k+1}^R(M,N),\;\;k\gg 0\,.\]
(It is enough to take $2k > n$, as resolutions in those degrees are $2$--periodic.)
This definition makes sense, as soon as the lengths appearing are finite. This
certainly happens if $R$ has an isolated singular point, but more generally
as soon as one of the two modules is locally free away from the singular point.\\

\begin{examples}
Hochster's Theta pairing is easy to compute. Let us take the simplest 
possible singularities, of type $A_1$ in two or three variables, and the simplest 
modules on them:

{\rm (i)} Take  $f=xy$, $M=\C\{x,y\}/(x)$, $N=\C\{x,y\}/(y)$, $K=\C\{x,y\}/(x,y)$. The
resolution of $M$ as $R=\C\{x,y\}/(xy)$-module is
\[\cdots \lra R \stackrel{x}{\lra}R \stackrel{y}{\lra}R \stackrel{x}{\lra}R \lra  M \lra 0\,,\]
thus, the ``matrix factorization'' is simply the factorization $x\cdot y$ of $f$ itself.
The $Tor_k^R(M,M)$ are computed by tensoring the above with $M$, that is, ``going mod $x$'', 
and we obtain
\[
\cdots \lra \C\{y\} \stackrel{0}{\lra} \C\{y\} \stackrel{y}{\lra}\C\{y\} \stackrel{0}{\lra} \C\{y\}  
\]
Hence we find 
\[length(Tor^{R}_{even}(M,M))=0,\;\;\;\;length(Tor^{R}_{odd}(M,M)) =1,\] 
so that $\theta(M,M)=-1$. To get $\theta(M,N)$, we have to ``go mod $y$'' instead
and obtain
\[\cdots \lra \C\{x\} \stackrel{x}{\lra} \C\{x\} \stackrel{0}{\lra}\C\{x\} \stackrel{x}{\lra} \C\{x\}  \]
so that now 
\[length(Tor^{R}_{even}(M,N))=1,\;\;\;\;length(Tor^{R}_{odd}(M,N))=0,\] hence 
$\theta(M,N)=1$. Similarly, by going ``mod $(x,y)$'', one finds each of
$Tor^{R}_{even/odd}(M,K)$ to be one-dimensional, whence $\theta(M,K)=0$.

{\rm (ii)} Take $f=xy-z^2$, $M=\C\{x,y,z\}/(x,z)$. (Note that $M=\calo_{L}$, where $L$ is the line
$x=z=0$ on the quadric cone $f=0$.) A  matrix factorisation $(A,B)$ associated to 
$M$ is given by
\[ A=\left(\begin{array}{cc} y&-z\\-z&x \end{array} \right),\;\;
B=\left(\begin{array}{cc} x&z\\z&y \end{array} \right)\,,\]
and $Tor^R_k(M,M)$ is the homology of the complex 
\[
\cdots \lra \C\{y\}^2 \stackrel{\alpha}{\lra} \C\{y\}^2 \stackrel{\beta}{\lra}\C\{y\}^2 
\stackrel{\alpha}{\lra} \C\{y\}\lra \cdots\,, \]
where
\[ \alpha=\left(\begin{array}{cc} y&0\\0&0 \end{array} \right),\;\;\beta=
\left(\begin{array}{cc} 0&0\\0&y \end{array} \right)\,.\]
So we find
\[length(Tor^{R}_{even}(M,M))=1,\;\;\;\;length(Tor^{R}_{odd}(M,M))=1,\] hence 
$\theta(M,M)=0$.\\
\end{examples}

We list some properties that follow immediately from the definition.\\

{\bf Properties of Hochster's Theta pairing:}\\
\begin{enumerate}[\rm(i)]
\item 
$\theta(M,N)=0$ if $M$ or $N$ is a free $R$-module.
\item
$\theta(M,N)=\theta(N,M)$.
\item
$\theta(M,N)$ is {\em additive over short exact sequences}:
if
\[ 0 \lra M' \lra M \lra M'' \lra 0\]
is a short exact sequence of $R$--modules, then
\[ \theta(M,N)=\theta(M',N) +\theta(M'',N)\,.\]
\item
$\theta(M,N)= -\theta(syz(M),N)$, where $syz(M)$ is the first
syzygy of $M$.
\item
$\theta(M,N)=0$ if either $M$ or $N$ is of finite projective dimension.
\item
$\theta(M,N)=0$ if either $M$ or $N$ is artinian.
\end{enumerate}

All properties, except maybe the last, are obvious. For (vi), note first that with
$m=(x_{0},...,x_{n})R$ the maximal ideal of the ring, 
\[\theta(M,R/m)=rank(F)-rank(G)=0\,,\] 
as tensoring a minimal resolution of $M$ with $R/m$ kills all differentials, and so
$length(Tor^{R}_{even}(M,R/m))=rank(F)$ and $length(Tor^{R}_{odd}(M,R/m))=rank(G)$.
For the general case, we can use induction on the length of $N$, as any artinian $N$ sits in 
an exact sequence
\[ 0 \lra R/m \lra N \lra N' \lra 0\,.\]

We see that $\theta$ descends to a symmetric pairing on the $K$-group
of the category $mod(R)$, divided out by the classes of the free and
artinian modules.\\

\begin{remark} In \cite{buchweitz} the quantity
\[h(M,N):=length(Ext_{R}^{even}(M,N))-length(Ext_{R}^{odd}(M,N))\]
was  studied and called the {\em Herbrand difference\/}, as the additive analogue of the 
{\em Herbrand quotient} that arises in group cohomology for representations of cyclic groups,
modules over the (integral) ``hypersurface'' $x^{n}=1$.
While it is more convenient for us here to work with $Tor$, one has the relation
\[ h(M,N)=\theta(M^*,N)\,,\]
where $M^*=Hom_R(M,R)$ is the dual module and $M$ is maximal Cohen-Macaulay, 
whence the two notions are indeed equivalent, taking into account that $M\cong M^{**}$.
\end{remark}

\section{The meaning of  $\theta(M,N)$}

An interesting case arises  when $M=\calo_Y=R/I; N=\calo_Z=R/J$,
where $Y, Z \subset X_0$ are subspaces of $X_0$, defined by ideals $I$ and $J$
respectively. By additivity over short exact sequences and using the fact that 
every module admits a finite filtration with sub-quotients of the form $ R/I$,
knowing all $\theta(\calo_Y,\calo_Z)$ determines $\theta(M,N)$ for all 
modules $M$ and $N$.\\

In analogy with Serre's Tor-formula, $\theta(\calo_Y,\calo_Z)$ should have something to 
do with intersections taking place inside $X$, the domain of $f$.  
The aim of this paper is to clarify this relation.

\begin{theorem}[{\cite{hochster}}]
In the above situation,
$$\theta(\calo_Y,\calo_Z)=I(Y,Z)$$
in case that $Y \cap Z =\{0\}$. Here $I(Y,Z)$ is the ordinary intersection
multiplicity of $Y$ and $Z$ {\em in the ambient smooth space} $(\C^{n+1},0)$.
\end{theorem}

{\bf Idea of proof:}
The result follows easily from Serre's Tor formula for the intersection multiplicity 
and the ``change of rings exact sequence''
\[
\xymatrix@R5pt{
\cdots\ar[r]&Tor^P_k(M,N) \ar[r]&Tor^R_k(M,N) \ar[r] &Tor^R_{k-2}(M,N) \ar[r]&\\
\ar[r]&Tor^P_{k-1}(M,N) \ar[r]&Tor^R_{k-1}(M,N) \ar[r] &Tor^R_{k-3}(M,N) \ar[r]&\cdots
}
\]
relating $Tor$ over $P$ and over $R$. \hfill $\Diamond$

\vskip 20pt
On the other hand, if $f \in P=\C[x_1,x_2,\ldots,x_{2m+2}]$ is a homogeneous
polynomial of degree $d$ in $2m+2$ variables, then $f$ defines a homogeneous 
cone
$X_0=f^{-1}(0) \subset {\C}^{2m+2}$ and an associated $2m$-dimensional
projective hypersurface
\[T:=V(f) \subset \P^{2m+1} \]
of degree $d$. If  $Y$ and $Z$ are homogeneous sub-varieties of $X_0$
of codimension $m$, then the projectivizations of $Y$ and $Z$ are codimension
$m$ cycles in $T$, whose fundamental classes in $H^m(T)$ we 
denote by $[Y]$ and $[Z]$, respectively. 
The graded version of Hochster's theorem just stated then yields the following result.\\

\begin{theorem} 
If $Y$ and $Z$ intersect transversely, then
\[  
\theta(\calo_Y,\calo_Z) = -\frac{1}{d}[[Y]]{\cdot}[[Z]]\,,
\]
where $[[Y]]:=d[Y]-\deg(Y){\cdot}h^m$ is the {\em primitive class} of $[Y]$, with $h\in H^{1}(T)$
the hyperplane class, and $\deg(Y)$ the degree of the subvariety $Y$ in $\P^{2m+1}$.\\
\end{theorem}

{\bf Proof:}
Recall that the primitive class of a cycle $Y$ is the projection of its fundamental class 
$[Y]\in H^{m}(T)$ into the orthogonal complement to $h^{m}$ with respect to the intersection 
pairing into $H^{2m}(T)\cong \C$. As $h^{m}h^{m} = d = \deg(T)$ and $[Y]{\cdot}h^{m}=\deg(Y)$, 
the description $[[Y]]=d[Y]-\deg(Y)h^{m}$ of the primitive class follows. Substituting, the claim 
can be reformulated as
\begin{align*}
\theta(\calo_Y,\calo_Z) &=  -\frac{1}{d}[[Y]]{\cdot}[[Z]] = -d[Y]{\cdot}[Z] + \deg(Y)\deg(Z)\,,
\end{align*}
where $[Y]{\cdot}[Z]$ denotes the intersection form on the cohomology of projective space.

The claim on Hochster's Theta pairing now follows from an argument on Hilbert--Poincar\'e series, 
the generating functions
$\HH(M) =\sum_{i}\dim (M_{i})t^{i}$, where $M=\oplus_{i\in\Z}M_{i}$ is a finitely generated 
graded module. For a complex of graded 
modules $C^{j}$ set  $\HH(C^{\bullet})=\sum_{j}(-1)^{j}\HH(C^{j})$.
The latter alternating sum is defined, as long as for a fixed degree $i$ 
only finitely many of the modules $C^{j}$ satisfy $C^{j}_{i}\neq 0$, and in that case 
$\HH(C^{\bullet})= \HH(H(C^{\bullet}))$.

If $C^{\bullet}=(\cdots C^{1}\to C^{0})\to M$ is a minimal homogeneous free resolution, then 
its terms $C^{j}$ are generated in higher and higher degrees as $j$ increases, whence 
$\HH(C^{\bullet})$ is summable. Exactness shows the alternating sum to equal $\HH(M)$.

If $N$ is a second finitely generated graded module and $D^{\bullet}\to N$ a 
minimal homogeneous free resolution, then the complex $C^{\bullet}\otimes_{R}D^{\bullet}$ that 
is the (total complex of the) tensor product of the resolutions still satisfies the summability 
condition. Its Hilbert--Poincar\'e series is independent of the choice of resolutions and 
denoted $\HH(M \stackrel{\L}{\tensor}_R N)$. It follows readily; see, for example, 
\cite[Lemma 7]{AvB}; that
\[
\HH(M \stackrel{\L}{\tensor}_R N)=\frac{\HH(M)\HH(N)}{\HH(R)}\,.
\]
Furthermore, in case that $M={\calo}_Y, N={\calo}_Z$, for cycles $Y,Z$ that intersect transversely,
$Tor_{i}^{R}(M,N)$ is of finite length for $i > 0$, and $Tor_{i+2}^{R}(M,N)$ is isomorphic
as graded vector space to $Tor_{i}^{R}(M,N)$, shifted in degrees by $d$. 
Thus, equating the Hilbert--Poincar\'e series of $M \stackrel{\L}{\tensor}_R N$ with that of its 
homology results in
\[
\HH(M{\stackrel{\mathbb L}{\tensor}_R} N)=\HH(M \tensor_{R} N) + 
polynomial+\frac{\HH(Tor_{ev}^{R}(M,N)){-}\HH(Tor_{odd}^{R}(M,N))}{1-t^d}\,.
\]
Next observe that $M\otimes_{R}N\cong M\otimes_{P}N$ and that $Tor_{i}^{P}(M,N)$ too is 
of finite length for $i > 0$, equal to zero for $i>2m+2$. Thus, $\HH(M \tensor_{R} N)$ and 
$\HH(M \stackrel{\L}{\tensor}_{P} N)$ differ only by a polynomial. Now compare the residues at $t=1$:\\
Because $Y,Z$ are of codimension $m$ in $X_{0}$ of dimension $2m+1$, 
the Hilbert--Poincar\'e series of $M,N$ are of the form
\begin{align*}
\HH(M) &= \frac{p_{M}}{(1-t)^{m+1}}\,,\quad \HH(N) = \frac{p_{N}}{(1-t)^{m+1}}\,,
\intertext{with $p_{M}, p_{N}\in\Z[t]$ satisfying $p_{M}(1) = \deg(Y), p_{N}(1) = \deg(Z)$, whence}
\HH(M \stackrel{\L}{\tensor}_R N)&= \frac{\HH(M)\HH(N)}{\HH(R)} = \frac{p_{M}p_{N}}{1-t^{d}}\,,
\end{align*}
as $\HH(R) =(1-t^{d})/(1-t)^{2m+2}$. Therefore, the residue evaluates to
\begin{align*}
res_{t=1}\HH(M \stackrel{\L}{\tensor}_R N) = \frac{1}{d}p_{M}(1)p_{N}(1) =\frac{1}{d}\deg(Y)\deg(Z)\,.
\end{align*}
For the Hilbert--Poincar\'e series $\HH(M \tensor_{R} N)$ we get
\begin{align*}
res_{t=1}\HH(M \tensor_{R} N) = res_{t=1}\HH(M \stackrel{\L}{\tensor}_{P} N) = [Y]\cdot[Z]\,.
\end{align*}
Finally,
\begin{align*}
res_{t=1}\frac{\HH(Tor_{ev}^{R}(M,N)){-}\HH(Tor_{odd}^{R}(M,N))}{1-t^d} 
&= \frac{1}{d}\theta(\calo_Y,\calo_Z)\,.
\end{align*}
Putting it all together, the equality of residues becomes
\[
\frac{1}{d}\deg(Y)\deg(Z) = [Y]\cdot[Z] +  \frac{1}{d}\theta(\calo_Y,\calo_Z)
\]
and solving for the Theta pairing yields the claim.
\hfill $\Diamond$

\begin{example}
Let $f\in \C[x_{1},x_{2},x_{3},x_{4}]$ be the equation of a cubic surface $S$ in $\P^{3}$.
A line $L$ on $S$ is given by two linear forms $l_{1},l_{2}$ such that
$f=l_{1}q_{1}+l_{2}q_{2}$, for suitable quadratic polynomials $q_{1},q_{2}$. 

The matrix factorization associated to $L$ has the form
\[ 
A=\left(\begin{array}{cc} l_{1}&-q_{2}\\l_{2}&q_{1}\end{array} \right)\;,\;
B=\left(\begin{array}{cc} q_{1}&q_{2}\\-l_{2}&l_{1} \end{array} \right)\,,
\]
and one easily determines from this $\theta(\calo_{L},\calo_{L'})$ for a pair of lines on $S$. 
The Theta pairing then recovers the $E_{6}$--lattice on the primitive cohomology 
of a smooth cubic surface from the configuration of the lines on $S$. 
In fact, for lines $L,L'$ one has the following table of dimensions of the torsion groups and 
values of Hochster's Theta pairing:
\[
\begin{array}{|c|c|c|c|}
\hline position&skew&transverse&identical\vphantom{\bigg(}\\
\hline\vphantom{\bigg(}
Tor^{R}_{even}(\calo_{L},\calo_{L'})&1&\hphantom{-}0&\hphantom{-}0\\
\vphantom{\bigg(}
Tor^{R}_{odd}(\calo_{L},\calo_{L'})&0&\hphantom{-}2&\hphantom{-}4\\
\hline
\vphantom{\bigg(}
\theta(\calo_{L},\calo_{L'})&1&-2&-4\\
\hline
\end{array}
\]
\smallskip

and in each case this agrees with the value predicted by the preceding result, even when the cycles
are not transversal.
\end{example}

\vskip 10pt
If, however, $f$ is {\em not} (quasi-) homogeneous, there no longer will be a
projective variety to do intersection theory on. So the question arises as to
the meaning of $\theta$ in the general case. For this, the geometry of
the link $L=X_0 \cap S$ of the isolated singularity inside the Milnor sphere will be relevant.
We will continue with the case $n=2m+2$. The non-vanishing cohomology groups of
$L$ then are:
$$\Z=H^0(L)  \;\;,\;\;\; H^{2m}(L)\;\;,\;\;\;  H^{2m+1}(L) \;\;,\;\;\; 
H^{4m+1}(L)=\Z\,.$$
The fundamental class of a codimension $m$ cycle on $X_{0}$ will land in $H^{2m}(L)$, 
but if we want to attach a number to two classes $[A],[B] \in H^{2m}(L)$ something
new has to be involved, because the dual space to $H^{2m}(L)$ is $H^{2m+1}(L)$.
J.~Steenbrink came up with the following conjecture:\\

\begin{conjecture}{\em (\cite{steenbrink}):\/}
$$\theta(\calo_A,\calo_B)=lk([A],[B])\,.$$

\end{conjecture}
Here $lk:H^{2m}(L) \times H^{2m}(L) \lra \Z$ is the co-called {\em linking 
form}, while $[A] \in H^{2m}(L)$ is the topological fundamental class obtained by intersecting
the cycle $A$ on $X_0$ with the Milnor sphere $S$.

The formula states that, geometrically, one takes the classes 
$[A]$ and $[B]$, then shifts $[A]$ transverse to $L$ into $S=S^{4m+3}$ to get 
$\widetilde{[A]}$. The cycles $\widetilde{[A]}$ and ${[B]}$ now have the right 
codimensions to link in $S$. To shift, the canonical trivialization of 
the normal bundle to $L$ determined by the values of the function $f$ is used.

\begin{example} Let us consider the earlier example, where $M=\C\{x,y\}/(x)$, $N=\C\{x,y\}/(y)$ 
on the $A_1$ singularity $xy=0$. The classes $[M]$ and $[N]$ are geometrically represented 
by circles $x=0, |y|=\epsilon$ and $y=0, |x|=\epsilon$ respectively,
with their standard anti-clockwise orientation. The linking number of these 
circles is $+1$. The cycle $\widetilde{[M]}$ is represented by the circle $xy=t$, for $|t|\neq 0$ 
small and fixed, $|y|$ fixed. We see that if $y$ runs anti-clockwise, then
$x$ runs clockwise, so the linking number of $\widetilde{[M]}$ with $[N]$ is $-1$.
This is in accordance with the computation of Hochster's Theta pairing.
\end{example}

Originally, J.~Steenbrink formulated the conjecture in terms of the 
{\em variation mapping in the Milnor fibration}, but this is equivalent to the 
above conjecture by standard results on the topology of isolated hypersurface singularities, 
see section 5. We note that the conjecture is compatible with 
the special cases covered by theorem 1 and 2. The following is the main result of this paper:\\

{\bf Main Theorem} Let $f \in P:=\C\{x_0,x_1,\ldots,x_n\}$ define an 
isolated singularity, and let $M$ and $N$ be $R=P/(f)$-modules.
\begin{enumerate}[\rm(i)]
\item If $n$ is odd, then $\theta(M,N)=0$. 
\item
If $n$ is even, then
\[ 
\theta(M,N)=lk(ch(M),ch(N))\,.
\]
\end{enumerate}
Here $ch: K^0(L) \lra H^{ev}(L)$ is the Chern-character.
Only the $2m$--component $ch^{2m} \in H^{2m}(L;\Q)$ contributes to the 
linking, so that alternatively we might write 
$\theta(M,N)=lk(ch^{2m}(M),ch^{2m}(N))$.

\section{Interpretation in Algebraic  K-theory}
\bigskip
We start with a description of $\theta(M,N)$ in terms of {\em algebraic} 
$K$-theory. After Quillen, K-theory can be defined for any abelian  
(or even exact) category $\cala$.  In either case,  $K_0(\cala)$ is just the 
Grothendieck group of $\cala$, so elements are represented as formal differences
$[X]-[Y]$ of isomorphism classes of objects in $\cala$ modulo relations 
$0=[X]-[Y]+[Z]$ for each exact sequence $0 \lra X \lra Y \lra Z \lra 0$.

A fundamental result obtained by Quillen \cite{quillen} is the {\em localization
sequence}: a Serre-subcategory $\cals$ of $\cala$ gives rise to a
long exact sequence of higher $K$-groups
\[ 
\cdots \lra K_i(\cals) \lra K_i(\cala) \lra K_i(\cala/\cals) \stackrel{\partial}{\lra} 
K_{i-1}(\cals) \lra K_{i-1}(\cala) \lra \cdots
\]

Elements in higher $K$-groups are harder to describe. 
Gillet-Grayson \cite{gilletgrayson}, identify the higher groups as 
$K_i(\cala)=\pi_i(G\cala)$, where $G \cala$ is a certain
simplicial space, whose $k$-simplices are given by pairs of flags of objects in $\cala$
\[ X_0 \inj X_1 \inj \cdots \inj X_k\]
\[ Y_0 \inj Y_1 \inj \cdots \inj Y_k\]
with compatible identifications $X_j/X_i \approx Y_j/Y_i$.
Nenashev \cite{nenashev} has shown, building on earlier work of Sherman \cite{sherman1}, 
that all elements in $K_1(\cala)$ can be represented by so-called {\em double short 
exact sequences} (d.s.e.s): these are pairs of exact sequences on the same three objects 
of $\cala$: 
\[
{\xymatrix@R.8pt{
&0 \ar[r]& A\ar[r]^-{f}\ar@{=}[dd]& B \ar[r]^-{g}\Equal[dd]&C \ar[r]\Equal[dd] & 0\\
 \Xi:=   &      \\
&0 \ar[r]& A\ar[r]^-{h}& B \ar[r]^-{k}&C \ar[r] & 0}}
\]
(Of course, the diagrams are {\em not} supposed to commute.)

\begin{examples} 
{\rm (i)} The elements where $A=0$ lead to special $K_1$-elements, associated to a pair of isomorphisms, or to an automorphism of a single object. In the set-up of
\cite{nenashev}, an automorphism $\beta: B \lra B$ corresponds to the d.s.e.s
\[
{\xymatrix@R2.5pt{
&0 \ar[r]& 0 \ar[r]\Equal[dd]& B \ar[r]^-{Id}\Equal[dd]&B \ar[r]\Equal[dd] & 0\\
\\
&0 \ar[r]& 0\ar[r]& B \ar[r]^-{\beta}_{\cong}&B \ar[r] & 0}}
\]

{\rm (ii)} By a {\em cyclic diagram} in $\cala$ we mean a diagram%
\footnote{Although $X=Y$, we prefer to distinguish these two copies of the same object to indicate
clearly in the following diagrams, where each copy comes from.}
 of the form
\[
\xymatrix@R.6pt{
&0 \ar[r]& A\ar[r]^-{\alpha}& X \ar[r]^-{\beta}\Equal[dd]&C \ar[r]\ar[dd]^-{\nu} & 0\\
 \xi:=   &      \\
&0 &\ar[l] D\ar[uu]^-{\mu}&\ar[l]_-{\delta} Y &\ar[l]_-{\gamma} B & \ar[l] 0
}
\]
where  both rows are short exact sequences in $\cala$. 
In case that $\mu$ and $\nu$ are {\em isomorphisms}, such a diagram determines
a class $\{\xi\} \in K_1(\cala)$.
As some maps are going in the wrong direction, we do not have
a d.s.e.s, but one can obtain one by considering the diagram

\[
\xymatrix
@R7pt
{
&0 \ar[r]& A\oplus B\ar[r]^-{p}& X\oplus B\oplus D \ar[r]^-{q}\Equal[dd]&
C\oplus D \ar[r]\ar[dd]^-{n}_{\cong} & 0\\
\\
&0 \ar[r] &D\oplus B\ar[uu]^-{m}_{\cong}\ar[r]^-{r} &Y \oplus B\oplus D\ar[r]^-{s}& B\oplus D \ar[r]  &0
}
\]
where the maps are in block-form, defined as follows

\[
p:=\left(\begin{array}{cc}
\alpha&0\\0&1\\0&0\\
\end{array}\right),\;\;\;q:=\left(\begin{array}{ccc}\beta&0&0\\0&0&1\\\end{array} \right),\;\;\; r:=\left(\begin{array}{cc}0&\gamma\\0&0\\1&0\\\end{array}\right),\;\;\;s:=\left(\begin{array}{ccc}0&1&0\\\delta&0&0\\\end{array}\right)
\]
and $m:=\left(\begin{array}{cc}\mu&0\\0&1\end{array}\right),\;\;n=\left(\begin{array}{cc}\nu&0\\0&1\end{array}\right)$. Thus we get a d.s.e.s 
\[
\xymatrix
@R7pt
{
&0 \ar[r]& A\oplus B\ar[r]^-{p}\Equal[dd]& X\oplus B\oplus D \ar[r]^-{q'}\Equal[dd]&
B\oplus D \ar[r]\Equal[dd] & 0\\
\Xi:=&\\
&0 \ar[r] &A\oplus B\ar[r]^-{r'} &Y \oplus B\oplus D\ar[r]^-{s}& B\oplus D \ar[r]  &0
}
\]
where $q':=n q, r':=r m^{-1}$. So it determines a class
\[ \{\xi \}:=[\Xi] \in K_1(\cala)\,.\]
{\rm (iii)} A $2$-periodic complex $C^{\bullet}$ in $\cala$:
\[ \cdots \stackrel{a}{\lra} Y \stackrel{b}{\lra} X \stackrel{a}{\lra} Y 
\stackrel{b}{\lra} X \stackrel{a}{\lra} \cdots\]
(where $Y$ is on the even spots and $X=Y$)
determines a canonical cyclic diagram
\[
\xymatrix@R.6pt{
&0 \ar[r]& A\ar[r]^-{\alpha}& X \ar[r]^-{\beta}\Equal[dd]&C \ar[r]\ar[dd]^-{\nu} & 0\\
 \xi:=   &      \\
&0 &\ar[l] D\ar[uu]^-{\mu}&\ar[l]_-{\delta} Y &\ar[l]_-{\gamma} B & \ar[l] 0
}
\]
where $A=Ker(a)$, $C=Im(a)$, $B=Ker(b)$, $D=Im(b)$.\\ 
If the complex is exact, then $\mu$ and $\nu$ are isomorphisms, and so
$C^{\bullet}$ determines the class
\[\{C^{\bullet}\} :=\{\xi\} \in K_1(\cala)\]
\end{examples}

{\bf The boundary map.} Sherman \cite{sherman2} has given an explicit description for the
boundary map
\[ K_1(\cala) \stackrel{\partial}{\lra} K_0(\cals)\]
that we will need. If $f: A\lra B$ is a morphism in $\cala$, let us
denote, as in \cite[p.177]{sherman2}, by an overline the corresponding map in
$\cala/\cals$: $\overline{f}:\overline{A}\lra \overline{B}$.\\
If $f: A \lra B$ in $\cala$ is an $\cals$-isomorphism, that is, $Ker(f), Coker(f) \in \cals$, we put
\[\chi(f):=[Coker(f)]-[Ker(f)] \in K_0(\cals)\]
If $A,B$ in $\cala$  are lifts of $\overline{A},\overline{B}$ in 
$\cala/\cals$ and $\phi:\overline{A}\lra \overline{B}$ is
an $\cala/\cals$-morphism, then it can be ``lifted'' to a diagram 
\[
\xymatrix{
A&\ar[l]_-{a} C\ar[r]^-{f}&D&\ar[l]_-{b}B\,,
}
\]
where $a$ and $b$ are $\cals$-isomorphisms and $\phi=\overline{b}^{-1}\overline{f}\overline{a}^{-1}$.
In case $f$ itself  is an $\cals$-isomorphism, one puts 
\[ \chi(\phi):=\chi(f)-\chi(a)-\chi(b) \in K_0(\cals)\,.\]
This class is independent of the choices made.

Given an element $x \in K_1(\cala/\cals)$ represented by
a double short exact sequence
\[
{\xymatrix@R2.5pt{
&0 \ar[r]& \overline{A}\ar[r]^-{\lambda}\ar@{=}[dd]& \overline{B} \ar[r]^-{\mu}\Equal[dd]&
\overline{C} \ar[r]\Equal[dd] & 0\\
\\
&0 \ar[r]& \overline{A}\ar[r]^-{\sigma}& \overline{B} \ar[r]^-{\tau}&\overline{C} \ar[r] & 0}}
\]
one can ``lift'' the morphisms involved in such a way%
\footnote{If $\lambda= \overline{b}^{-1}\overline{\ell}(\overline{a}^{-1})$ and 
$\mu=\overline{c}^{-1}\overline{g}\overline{b'}^{-1}$, take $\alpha_{1}= \overline{a}, 
\beta_{1}=\overline{b}^{-1}$ and
$\overline{m}=\overline{g}\overline{b'}(\overline{b}^{-1}), \gamma_{1}=\overline{c}^{-1}$ and so on.}
as to obtain an extended diagram in 
$\cala/\cals$ with exact rows and commuting squares at the top and bottom,
\[
\xymatrix@R7pt{
&0 \ar[r]& \overline{A_{1}}\ar[r]^-{ \overline{\ell}}\ar[dd]_-{\alpha_{1}}^-{\cong}&
\overline{B_{1}} \ar[r]^-{ \overline{m}}\ar[dd]^{\beta_{1}}_-{\cong}&
\overline{C_{1}} \ar[r]\ar[dd]^{\gamma_{1}}_-{\cong} & 0\\
\\
&0 \ar[r]& \overline{A}\ar[r]^-{\lambda}\ar@{=}[dd]& \overline{B} \ar[r]^-{\mu}\Equal[dd]&
\overline{C} \ar[r]\Equal[dd] & 0\\
\\
&0 \ar[r]& \overline{A}\ar[r]^-{\sigma}& \overline{B} \ar[r]^-{\tau}&\overline{C} \ar[r] & 0\\
\\
&0 \ar[r]& \overline{A_{2}}\ar[r]^-{ \overline{s}}\ar[uu]^-{\alpha_{2}}_-{\cong}&
\overline{B_{2}} \ar[r]^-{ \overline{t}}\ar[uu]_-{\beta_{2}}^-{\cong}&
\overline{C_{2}} \ar[r]\ar[uu]_-{\gamma_{2}}^-{\cong} & 0}
\]

\begin{theorem} {\rm (Theorem 2.3 of \cite{sherman2})}\\
The boundary map $K_1(\cala/\cals) \stackrel{\partial}{\lra} K_0(\cals)$ sends the 
element $x$ to
\[ \partial(x)=\chi(\alpha)-\chi(\beta)+\chi(\gamma)  \in K_0(\cals)\,,\] 
where
\[\alpha=\alpha_2^{-1}\alpha_1,\;\;\;\beta=\beta_2^{-1}\beta_1,\;\;\gamma=\gamma_2^{-1}\gamma_1\,.\]
\end{theorem}

\begin{corollary} If in a cyclic diagram
\[
\xymatrix@R7.5pt{
&0 \ar[r]& A\ar[r]^-{\alpha}& X \ar[r]^-{\beta}\Equal[dd]&C \ar[r]\ar[dd]^-{\nu} & 0\\
\\
&0 &\ar[l] D\ar[uu]^-{\mu}&\ar[l]_-{\delta} Y &\ar[l]_-{\gamma} B & \ar[l] 0
}
\]
in $\cala$ the morphisms $\mu,\nu$ are $\cals$-isomorphisms, then the reduction
\[
\xymatrix@R.6pt{
&0 \ar[r]& \overline{A} \ar[r]^-{\overline{\alpha}}& 
\overline{X} \ar[r]^-{\overline{\beta}}\Equal[dd]&
\overline{C} \ar[r]\ar[dd]^-{\overline{\nu}}_-{\cong} & 0\\
 \xi:=   &      \\
&0 &\ar[l] \overline{D}\ar[uu]^-{\overline{\mu}}_{\cong}&
\ar[l]_-{\overline{\delta}} \overline{Y} &\ar[l]_-{\overline{\gamma}} 
\overline{B} & \ar[l] 0
}
\]
in $\cala/\cals$ satisfies: 
\[
\partial(\{\xi\})=([Coker(\nu)]-[Ker(\nu)])-([Coker(\mu)]-[Ker(\mu)])\,.\]
\end{corollary}

{\bf Proof:} The d.s.e.s associated to $\xi$ comes with given
lifts to $\cala$ from the cyclic diagram in $\cala$. The result now
follows by a direct application of Sherman's theorem to the associated double short
exact sequence. \hfill $\Diamond$ 
\vskip 10pt
{\bf Construction for hypersurfaces.}
We will apply the above theory to our situation of modules on
a hypersurface ring $R$ with an isolated singular point. We
let  $\cala:=mod(R)$, the category of finitely generated $R$-modules
and $\cals=art(R)$, the Serre-subcategory of artinian $R$-modules, 
or what is the same, of those that are  supported at the singular point.
The category $\cala/\cals$ is the category obtained by ``localizing
away'' from the singular point and can be identified with the category
$coh(U)$ of coherent sheaves on the punctured spectrum 
$U:=X_0\setminus \{0\}$. We will write ${\bf K}^i(U)$ for $K_i(coh(U))$.
As $U$ is smooth, it is the same as $K_i(Vect(U))$, where $Vect(U)$ is
the category of locally free $\calo_U$-modules.\\

{\bf Definition:} For a maximal Cohen-Macaulay module $M$ on a hypersurface
ring $R$ we put:
\[ [M]=[M_U] \in {\bf K}^0(U)\]
\[ \{M\}:=\{C^{\bullet}(M)_U\} \in {\bf K}^1(U)\]

where $(.)_U$ denotes restriction to $U$.\\

We recall that there exists a natural product
\[ {\bf K}^1(U) \times {\bf K}^0(U) \lra {\bf K}^1(U)\,,\]
induced in the obvious way by the tensor product $\otimes$ of modules,
and a trace map 
\[ {\chi}:=\ell \circ \partial: {\bf K}^1(U) \lra \Z \]
that is the composition of the boundary map 
${\bf K}^1(U) \stackrel{\partial}{\lra} {K}_0(art(R))$ and the isomorphism 
$\ell:K_0(art(R)) \lra \Z,\;\;\; M \mapsto length(M)$.\\

The following crucial result expresses Hochster's Theta pairing in terms of 
algebraic $K$-theory.

\begin{theorem}
\[\theta(M,N)={\bf \chi}(\{M\} \otimes [N])\,.\]
\end{theorem}

{\bf Proof:} Recall the exact $2$-periodic complex of $R$-modules
\[
\xymatrix{
C^{\bullet}(M):&\cdots\ar[r]^-{B}&G\ar[r]^-{A}&F\ar[r]^-{B}&G\ar[r]^-{A}&F\ar[r]^-{B}&\cdots
}
\]
that is defined by (a matrix factorisation $(A,B)$ attached to) $M$. 
By definition, the groups $Tor_{even/odd}^{R}(M,N)$
are obtained by tensoring this complex with $N$ and then taking homology. To put this differently,
consider the following cyclic diagram
\[
\xymatrix@R7.5pt{
&0 \ar[r]& Ker(A \otimes_{R} Id)\ar[r]& G \otimes_{R} N \ar[r]\Equal[dd]&
Im(A \otimes_{R} Id) \ar[r]\ar[dd]^-{\nu} & 0\\
\\
&0 &\ar[l] Im(B \otimes_{R} Id)\ar[uu]^-{\mu}&\ar[l] F\otimes_{R} N &
\ar[l] Ker(B \otimes_{R} Id)  & \ar[l] 0
}
\]
where $Id:N \lra N$ is the identity map. We then have
\[Tor_{even}^R(M,N)=Coker(\nu),\;\; Tor_{odd}^R(M,N)=Coker(\mu)\,.\]

On the other hand, the class $\{M\} \in {\bf K}^1(U)$ is represented by the cyclic diagram and 
resulting d.s.e.s. associated to the exact $2$-periodic complex
\[
\xymatrix{
\overline{C^{\bullet}(M)}:&\cdots\ar[r]^-{\overline{B}}&
\overline{G}\ar[r]^-{\overline{A}}&\overline{F}\ar[r]^-{\overline{B}}&
\overline{G}\ar[r]^-{\overline{A}}&\overline{F}\ar[r]^-{\overline{B}}&\cdots
}
\]
where the overline indicates the restriction to $U$. As we are now outside the singular locus,
the tensor product 
\[
\xymatrix{
\cdots\ar[r]&\overline{G} \otimes \overline{N}\ar[r]^-{\overline{A} \otimes \overline{Id} }&
\overline{F} \otimes \overline{N}\ar[r]^-{\overline{B} \otimes \overline{Id}}&
\overline{G} \otimes \overline{N}\ar[r]^-{\overline{A} \otimes \overline{Id}}&
\overline{F} \otimes \overline{N}\ar[r]^-{\overline{B} \otimes \overline{Id}}&\cdots
}
\]
of $\overline{C^{\bullet}(M)}$ with $\overline{N}$ on $U$ stays exact, and the class
$\{M\} \otimes [N] \in {\bf K}^1(U)$ is represented by the
cyclic diagram
\[
\xymatrix@R7.5pt{
&0 \ar[r]& Ker(\overline{A} \otimes \overline{Id})\ar[r]& 
\overline{G} \otimes \overline{N} \ar[r]\Equal[dd]&
Im(\overline{A} \otimes \overline{Id}) \ar[r]\ar[dd]^-{\overline{\nu}}_-{\cong} & 0\\
\\
&0 &\ar[l] Im(\overline{B} \otimes \overline{Id})\ar[uu]^-{\overline{\mu}}_-{\cong}&
\ar[l] \overline{F}\otimes \overline{N} &
\ar[l] Ker(\overline{B} \otimes \overline{Id})  & \ar[l] 0
}
\]
associated to that tensor product on $U$. 

It follows now from the corollary to the theorem of Sherman that
\[
\begin{array}{ccl}
\partial(\{M\}{\otimes} [N])&=&[Coker(\nu)]-[Coker(\mu)]\\
&=&[Tor^{R}_{even}(M,N)]-[Tor^{R}_{odd}(M,N)]
\end{array}
\]
and applying the trace map completes the argument. \hfill $\Diamond$

\section{Interpretation in terms of topological K-theory }

Topological $K$-theory is a generalized cohomology theory obtained from
topological $\C$-vector bundles. We refer to Atiyah's classic \cite{atiyah} for
a nice introduction. This theory fits into the above framework, if we take
for any topological space $X$ the abelian category $Vect_{top}(X)$ of topological $\C$-vector bundles on it: 
\[K^i(X)=K_i(Vect_{top}(X))\,.\]
As vector bundles can be pulled back, this is naturally a contravariant functor.
Distinctive features of the theory are:

$\bullet$ In $Vect_{top}(X)$ all (short) exact sequences can be split: if
$E \stackrel{i}{\inj} F \surj G$ is such a sequence, just pick a {\em hermitian metric} on $F$ 
and use orthogonal projection to obtain $F \cong E \oplus E^{\perp}, E^{\perp} \cong G$.

$\bullet$ One has {\em Bott periodicity}\,: $K^i(X) \cong K^{i+2}(X)$. So
we really have to consider only the two groups $K^0$ and $K^1$. 

$\bullet$ Elements
of $K^1$ can all be represented by {\em automorphisms} $E \stackrel{\alpha}{\lra} E$ of a 
vector bundle on $X$. If $x=[E \stackrel{\alpha}{\lra}E]$, then
$-x=[E \stackrel{\alpha^{-1}}{\lra} E]$.

$\bullet$ The elements in $K^{1}$ are {\em homotopy invariants}\,: if $t\mapsto \alpha_{t}$
is a path within the space of automorphisms of $E$, then
$x_{t}=[E \stackrel{\alpha_{t}}{\lra} E]$ is independent of $t$.

$\bullet$ There are  {\em Chern-class maps}
\[ch^{ev}: K^0(X) \lra H^{ev}(X)=\bigoplus_k H^{2k}(X,\Q)\]
\[ch^{odd}: K^1(X) \lra H^{odd}(X)=\bigoplus_k H^{2k+1}(X,\Q)\]
that are isomorphisms after tensoring with $\Q$.

$\bullet$ There is a version of Riemann-Roch in the differentiable context, stating that for a proper map, direct image in K-theory (defined via duality),
commutes with taking Chern-classes, up to multiplication with the Todd-class of the virtual tangent bundle of the map.\\ 

Of importance for us will be the special case of the constant map
$p: L \lra point$, where $L$ is the link, a compact, odd-dimensional manifold
with stably trivial tangent bundle. The induced map
\[  \chi_{top}:=p_*: K^1(L) \lra K^0(point)=\Z\] 
is the trace map in topological $K$-theory.\\
As $U$ and $L$ are homotopy equivalent, we will identify $K^i(U)$ with $K^i(L)$
without further mention, and write also $\chi_{top}:K^1(U) \lra \Z$, and so on.\\

We can compare algebraic and topological $K$-theory for $U$ using the obvious {\em topologification functor}
\[ top: Vect(U) \lra Vect_{top}(U)\]
sending  a locally free sheaf to its associated topological vector
bundle. This induces natural maps from algebraic to topological $K$--theory,
\[ 
\begin{array}{rcllccc}
{\mathbf K}^0(U) &\lra& K^0(U)=K^0(L)&\,, &[M]& \mapsto &[M]_{top}\,,\vphantom{\bigg(_{|}}\\
{\mathbf K}^1(U) &\lra& K^1(U)=K^1(L)&\,, &\{M\}& \mapsto &\{M\}_{top}\,.
\end{array}
\]\\

\begin{proposition}
\[ \theta(M,N)=\chi_{top}(\{M\}_{top} \otimes [N]_{top} )\]
\end{proposition}
{\bf Proof:} This follows from the naturality of topologification:
$(\{M\} \otimes [N])_{top} = \{M\}_{top}\otimes [N]_{top}$, and
compatibility of $\chi$ and $\chi_{top}$, that is, the commutativity
of
\[
\xymatrix@R18pt{
{\bf K}^1(U) \ar[r]^-{\chi}\ar[d]_-{top} &\Z \Equal[d]\\
K^1(U)\ar[r]^-{\chi_{top}}& \Z
}
\]
This ``well-known'' fact can be shown as follows. 
Take a finite map $p:X_0 \lra Z$, where $Z$ is smooth and contractible. 
We let $Z^*=Z \setminus \{0\}$, so that $p:U \lra Z^*$ is a finite (ramified)
covering. By functoriality of both ${\bf K^1}$ and ${K^1}$ and the fact that the
map to the point factors over $p$, we reduce to checking the statement for 
$Z^*$. The structure sheaf of the point $\{0\}$ is resolved by the Koszul 
complex; its class in either ${\bf K}^1(U)$ or ${K^1(U)}$ is mapped by 
$\chi$, respectively $\chi_{top}$, to $1 \in \Z$.
\hfill $\Diamond$\\

{\bf A different description of $\{M\}_{top}\ $.}
If $(A,B)$ is a matrix factorisation for a maximal Cohen-Macaulay module
$M$ on $X_{0}$, we can choose our good representative for $f: X\lra D$ in such a way
that the matrices $A$ and $B$ are holomorphic on $X$. Hence we have an
exact sequence of sheaves on $X$: 
\[
\xymatrix{
0\ar[r]&\calo_{X}^{p}\ar[r]^-{A}&\calo_X^p\ar[r]&M\ar[r] &0\,.
}
\]
(We will not make a notational distinction between objects over
the local ring and corresponding sheaves on $X$). As $\det(A)$ vanishes
only on the hypersurface $X_0$, the matrix $A(x)$ is an isomorphism for 
each $x  \in X^{*}=X\setminus X_{0}$. It determines hence a class
\[ \alpha(M) :=[\calo_{X^*}^p \stackrel{A}{\lra} \calo_{X^*}^p] \in K^1(X^*)\]
and, as $X^*$ retracts to $S-L$, we can also see
$\alpha(M) \in K^1(S-L)=K^1(X^*)$.
For a Milnor fibre $X_t$, we have inclusion maps $\partial X_t \inj X_t \inj X^*$, 
hence we get corresponding restriction maps
\[ K^1(X^*) \lra K^1(X_t) \lra K^1(\partial X_t)\] 
that send $\alpha(M)$ to  $\alpha(M)|_{X_t}$ and 
$\alpha(M)|_{\partial X_t}$ respectively. 
Furthermore, the local triviality of the Milnor fibration near the boundary
provides an identification $\rho_t: L \lra \partial X_t$, in particular, an
isomorphism $\rho_{t}^{*}:K^1(\partial X_t)\lra K^1(L)$.\\

\begin{theorem} For any $t \in D^*$, we have

\[\{M\}_{top}=\rho_t^*(\alpha(M)|_{\partial X_t})\,.\]
\end{theorem}
\bigskip
The proof will make use of the following lemma.\\

\begin{lemma} Choose a hermitian metric $(-,-)$ on $\calo_X^p$ and
let $B^{\dagger}$ be the hermitian adjoint of $B$. For $s \in \C^*$,
consider the matrix
\[A_s:=A-sB^{\dagger}\,.\]

{\rm (i)} The matrix $A_s(x)$ is invertible for all $x \in X_t$ if $t \notin \R_+{\cdot} s$, 
so defines a class $\alpha_s(M)|_{X_t} \in K^1(X_t)$, and that class satisfies 
\[ \alpha_s(M)|_{X_t}=\alpha(M)|_{X_t}\,.\]

{\rm (ii)} For each $s \neq 0$, the matrix $A_s(x)$ is invertible on $U$, so defines 
a class $\alpha_s(M)_{L} \in K^1(L)$. That class satisfies
\[ 
\alpha_s(M)_{L}=\rho_t^*(\alpha(M)|_{\partial X_t}) 
\] 
for $t \notin \R_+ {\cdot} s$.
\end{lemma}

{\bf Proof:} 
(i) Set $V_s:=\{  x \in X | \det(A_s(x))=0\}$. 
If $x \in V_s$, then there  exists by definition a vector 
$v=v(x) \neq 0$ in the kernel of $A_s(x)$. Hence, we have:
\[
A(x).v=s.B^{\dagger}(x).v\,.
\]
When we multiply this equation from the left with $B(x)$, we get 
\[s^{-1}f(x).v=B(x).B^{\dagger}(x).v\,,
\] 
that is, $s^{-1}f(x)$ is an {\em eigenvalue} of the matrix $B(x).B^{\dagger}(x)$. 
However, if $B(x).B^{\dagger}(x)v=\lambda.v$, then 
$(B^{\dagger}(x)v,B^{\dagger}(x)v)=\lambda.(v,v)$, so $\lambda \ge 0$. 
It follows that the image under $f$ of $V_s$ is contained in the half-line $\R_{+}{\cdot}s$ 
and consequently $V_{s}$ is disjoint from the Milnor fibre $X_t$ if $t \notin \R_+{\cdot} s$. 
Note that $X_t$ is then also disjoint from $V_{s'}$ for all $s' \in \R^{+}s$. 
So $s' \mapsto A_{s'}$ provides a continuous path from $A=A_0$ to $A_s$ inside the 
space of invertible matrices on the Milnor fibre $X_t$. Hence, $A$ and $A_s$ represent 
the same element in $K^1(X_t)$.

(ii) If $x \neq 0$, then the eigenvalue $\lambda=0$ cannot occur for
$B(x).B^{\dagger}(x)$: first, if $\lambda=0$ then $B^{\dagger}(x)v=0$, and hence also
$A(x).v=0$, as $v$ is in the kernel of $A_{s}(x)$. 
But then $0=(B^{\dagger}(x)v,w)=(v,B(x)w)$ shows that $v$ is 
orthogonal to $Im(B(x))= Ker(A(x))$, the last equality due to $x\in U$. 
As $v \in Ker(A(x))$, this shows $v=0$. 
We conclude that $A_s(x)$ is {invertible} for all $x \in U$ and thus
defines by restriction a class $\alpha_s(M)_L$. By construction, 
$\alpha_s(M)_L=\rho_t^*(\alpha_s(M)|_{\partial X_t})$, so the last statement follows
from (i).  \hfill $\Diamond$\\

{\bf Proof of the Theorem:}

We look at the exact 2-periodic complex $C^{\bullet}(M)_U$ {\em on U}:
\[  \cdots \stackrel{B}{\lra}  G \stackrel{A}{\lra} F \stackrel{B}{\lra} G \stackrel{A}{\lra} \cdots\]
where $F=G =\calo^p$ is the trivial rank $p$ bundle on $U$ and $F$ sits on the 
even spots. In order to keep the notation simple, we will omit the overline or 
notation $(-)_U$ to denote restriction to $U$; everything here takes place in $Vect_{top}(U)$.
We choose a hermitian metric on $F$ and split the exact sequences
\[ Ker(A) \inj G \surj G/Ker(A)\,,\quad Ker(B) \inj F \surj F/Ker(B)\,,\] 
to obtain isomorphisms  
\[G \cong Ker(A) \oplus G/Ker(A)\,,\quad F \cong F/Ker(B) \oplus Ker(B),\]
that then give rise to a diagram
\[
\xymatrix{
G\Equal[d]&\cong&Ker(A)&\oplus&G/Ker(A)\ar[d]^-{\displaystyle A}\\
F&\cong &F/Ker(B)\ar[u]^-{\displaystyle B}&\oplus &Ker(B)
}
\]
The class $\{M\}_{top}$ is represented by the
automorphism $G\lra F$ obtained from this diagram by inverting the
map $B$. Now, in general, if $E \stackrel{\alpha}{\lra} F$ is an isomorphism 
of vector bundles, and if $F \stackrel{\alpha^{\dagger}}{\lra} E$ is its
hermitian adjoint, then $\alpha^{\dagger}$ is {\em homotopic} to 
$\alpha^{-1}$, as the composition $\alpha \alpha^{\dagger}$ is positive definite, 
hence homotopic to the identity.
Accordingly our class $\{M\}_{top}$ is represented by the automorphism
$A+B^{\dagger}$. But this is precisely $\alpha_{-1}(M)_L$, and by part (ii) 
of the lemma it represents the same class as 
$\rho_t^*(\alpha(M)|_{\partial X_t})$.  \hfill $\Diamond$\\ 

\begin{corollary} If $n$ is {\em even\/} then $\theta(M,N)=0$ for all $M$ and $N$.
\end{corollary}
{\bf Proof:} 
If $n$ is even, then $K^1(X_t)=0$, as $X_t$ has the homotopy type of
a wedge of even dimensional spheres. Hence, the class of $M$ trivially satisfies
$\alpha(M)|_{X_t}=0$ in $K^1(X_t)$, and so certainly its restriction to the boundary vanishes
also: $\alpha(M)|_{\partial X_t}=0$.
But then 
\[\{M\}_{top}=\alpha_s(M)_{L}=\rho_t^*(\alpha(M)|_{\partial X_t})=0\,,\]
so $\theta(M,N)=\chi_{top}(\{M\}_{top}\otimes [N]_{top})=0$\,! \hfill $\Diamond$  

\section{The Linking Form}
We review here the basic properties of the {\em linking pairing} on the (co-)hom\-ol\-ogy of the link of an isolated hypersurface singularity.
More details can be found in \cite{agv}, \cite{levine}.\\
Given two disjoint $n$-dimensional cycles $\alpha,\beta$ in the Milnor sphere
$S=S^{2n+1}$, we can form the {\em linking number} $\ell(\alpha,\beta) \in \Z$,
which is defined as the intersection number $\Gamma \cdot \beta$ between 
a chain $\Gamma$ with $\partial \Gamma=\alpha$ and $\beta$. 
One has $\ell(\alpha,\beta)=(-1)^{n+1}\ell(\beta,\alpha)$, so that linking is symmetric for 
odd dimensional cycles. Consider the Milnor fibration 
$f: X^* \lra D^*$ as before. Fix $t \in D^*$ and 
use parallel transport along an anti-clockwise half-turn from $t$ to $-t$ to 
define a ``half-monodromy map'' 
\[h_{1/2}:H_n(X_t) \lra H_n(X_{-t})\,.\]
If $\alpha, \beta \in H_n(X_t)$, then the 
cycles  $\alpha$ and $h_{1/2}(\beta)$ are disjoint and have the
appropriate dimension to link in the $2n+1$ dimensional Milnor
sphere $S$. The resulting {\em Seifert form} of the singularity is
\[S: H_n(X_t) \times H_n(X_t) \lra \Z,\;\;(\alpha, \beta) \mapsto \ell(\alpha, h_{1/2}(\beta))\]  
If we restrict the Seifert form $S$ to $H_n(L)\cong H_n(\partial X_t) \subset H_n(X_t)$ we 
obtain a $(-1)^{n+1}$-symmetric form
\[lk: H_n(L) \times H_n(L) \lra \Z\,,\]
that we call the {\em linking form of the link}. So, geometrically, 
$lk(\alpha,\beta)=\ell(\tilde{\alpha}, \tilde{\beta})$, where $\tilde{\alpha}$ 
and $\tilde{\beta}$ are obtained from $\alpha$ and $\beta$
using the identification $H_n(L) \approx H_n(\partial X_{t})$ and
$H_n(L) \approx H_n(\partial X_{-t})$ respectively. Alternatively, we may
say that $\tilde{\alpha}$ and $\tilde{\beta}$ are obtained by ``pushing-aside''
$\alpha$ and $\beta$ in opposite directions, using the trivialisation of
the normal bundle of $L$ defined by $f$. It is clear that in fact one
only needs to push aside one of the cycles, so that $lk(\alpha,\beta)=\ell(\tilde{\alpha},\beta)$.

The geometric monodromy can be taken to be the identity on $\partial X_t$, so 
one can define a {\em variation mapping} 
\[Var: H_n(X_t,\partial X_t) \lra H_n(X_t),\]
obtained by mapping a relative cycle $\gamma$ to $[\gamma-h(\gamma)]$, where
$h$ is the monodromy. It is related to the Seifert form and the intersection 
pairing 
\[ (-,-) : H_n(X_t,\partial X_t) \times H_n(X_t) \lra \Z\]
by the formula $S(Var(\alpha),\beta)=(\alpha,\beta)$, and as the variation
is an isomorphism, we can write $S(\alpha,\beta)=(Var^{-1}\alpha,\beta)$.\\

{\bf Formulation in cohomology and K-theory.}
We will reformulate the above procedure in cohomological terms and cover it
through  a description in topological $K$-theory. We will restrict to the case
$n$ odd, so that the form $lk$ is symmetric. We write $n+1=2m+2$, so that
the non-vanishing cohomology of $L$ sits in degrees $0, 2m, 2m+1,4m+1$, and
$S$ has dimension $4m+3$. Written in cohomology, the linking form becomes a pairing on
the cohomology in degree $2m$:
\[ lk: H^{2m}(L) \times H^{2m}(L) \lra \Z\,.\]
From the long exact cohomology sequence of the pair $(S,S-L)$ and using the
fact that $S$ is a sphere, we get that the {coboundary} map $\delta$ 
is an {\em isomorphism}: 
\[ \delta:H^{2m+1}(S-L) \stackrel{\simeq}{\lra} H^{2m+2}(S,S-L)\] 
Furthermore, we have the {\em Thom isomorphism}
\[ t: H^{2m+2}(S,S-L) \stackrel{\simeq}{\lra} H^{2m}(L)\]
Combining the two, we get the {\em Alexander duality} isomorphism
\[ \lambda: H^{2m+1}(S - L) \stackrel{\simeq}{\lra}  H^{2m}(L)\]
Finally, we have a ``push aside map'' $\rho:=\rho_t: L \lra \partial X_t \subset S \setminus L$
that induces a map $\rho^*:H^{2m+1}(S-L) \lra H^{2m+1}(L)$. 
Combined we obtain a map
\[\gamma:=\rho^*\circ \lambda^{-1}: H^{2m}(L) \lra H^{2m+1}(L)\,.\]
Let finally
\[ <-,->: H^{2m+1}(L) \times H^{2m}(L) \lra \Z\]
be the Poincar\'e pairing of the oriented $4m+1$-manifold $L$. We then
have the following formula. 

\begin{proposition}
\[ lk(\alpha,\beta)=<\gamma(\alpha),\beta>\,.\]
$\hfill \Diamond$
\end{proposition}

We mimic  the above construction in topological K-theory.
The pair $(S,S-L)$ defines a long exact sequence in $K$-theory as it does in cohomology.
Furthermore, we have a {\em Thom-isomorphism in $K$-theory}, 
$T: K^0(S,S - L) \lra K^0(L)$, and the boundary map 
$\Delta:K^1(S-L) \to K^{0}(S,S-L)$ is an isomorphism as $S$ is an 
odd-dimensional sphere.  Comparing with (rational) cohomology, we 
obtain a diagram with commuting squares:
\[
\xymatrix{
K^1(S-L) \ar[r]^-{\Delta}_-{\simeq} \ar[d]_-{ch^{2m+1}}& K^0(S,S-L) \ar[r]^-{T}_-{\simeq}\ar[d]& 
K^0(L) \ar[d]^-{ch^{2m}}\\  
H^{2m+1}(S-L)\ar[r]^-{\delta}_-{\simeq} & H^{2m+2}(S,S-L) \ar[r]^-{t}_-{\simeq}& H^{2m}(L)}
\]

\begin{proposition} For a maximal Cohen-Macaulay module $M$ with classes
\[[M]_{top} \in K^0(L)\,,\quad \{M\}_{top} \in K^1(L)\,,\]
one has
\[ch^{2m+1}(\{M\}_{top})=\gamma(ch^{2m}([M]_{top}))\,.\]
\end{proposition}
{\bf Proof:} The Thom isomorphism in $K$-theory maps 
$E \stackrel{\alpha}{\lra} F$ on $S$, which is an isomorphism 
over $S-L$, to the index-bundle $[Coker(\alpha)]-[Ker(\alpha)]$.
So if we start with the $\calo_X$--resolution $0 \lra \calo^p_X \lra \calo^p_X \lra M \lra 0$ of 
$M$, then $[M]_{top} \in K^0(L)$
is just the image of $\alpha(M) \in K^1(S-L)$ under $T\circ \Delta$.
The geometric map $\rho_t: L \lra \partial X_t$ induces a 
commutative diagram
\[
\xymatrix{
K^1(S-L) \ar[r] \ar[d]_-{ch^{2m+1}}& K^1(\partial X_t) \ar[r]^-{\rho_t^*}& 
K^1(L)\ar[d]^-{ch^{2m+1}}\\
H^{2m+1}(S-L) \ar[r] & H^{2m+1}(\partial X_t) \ar[r]^-{\rho_t^*}& H^{2m+1}(L)
}
\]
By Theorem 4.2, the element $\alpha(M) \in K^1(S-L)$ gets mapped to $\{M\}_{top} \in K^1(L)$. The 
formula now follows from the definition of 
$\gamma=\rho_t^* \circ \delta^{-1} \circ t^{-1}$.\hfill $\Diamond$

\begin{corollary} (Proof of the main result)
\[
\begin{array}{cccl}
\theta(M,N)&=&\chi(\{M\} \otimes [N])&\textup{by Theorem 3.4\,,}\\
        &=&\chi_{top}(\{M\}_{top} \otimes [N]_{top})&\textup{by Proposition 4.1\,,}\\
        &=&<ch(\{M\}_{top}),ch([N]_{top})>&\\
        &=&<\gamma(ch([M]_{top})),ch([N]_{top})>&\textup{by Proposition 5.2\,,}\\
        &=&lk(ch([M]_{top}),ch([N]_{top}))&\textup{by Proposition 5.1.}\\
\end{array}
\]
$\hfill\Diamond$
\end{corollary}
\begin{remark}
{\em
Given a matrix factorisation $(A,B)$ for a maximal Cohen-Mac\-au\-lay module
$M$ one can find de Rham representatives for the Chern-classes that lie in 
the {\em unstable cyclic homology} of the hypersurface.\\ 
Consider $P=\C\{x_0,x_1,\ldots,x_n\}$ as a module over $\C\{t\}$ with
$t$ acting as multiplication by $f \in P$.
Denote by $\Omega^p$ the module of germs of $p$-forms on $\C^{n+1}$ and
let $\Omega^p_f:=\Omega^p/(df\wedge \Omega^{p-1})$ be the module of 
relative differentials. 
One puts
\[\omega(M):=dA \wedge dB\]
The components of the chern-character 
\[ch_M:=tr(exp(\omega(M)))=\sum_i \frac{1}{i!}\omega^i(M)\]
are well-defined classes 
\[\omega^i(M) := tr((dA \wedge dB)^{i}) \in \Omega^{2i}/(df\wedge d\Omega^{2i-1})\,.\]
There are, however, also {\em odd-degree\/} classes
\[\eta^i(M):=tr(A dB (dA \wedge dB)^i) \in \Omega^{2i+1}_f/d\Omega_f^{2i}\,.\]
The group $\Omega^{2i+1}_f/d\Omega_f^{2i}$ can be identified with the 
{\em cyclic homology} $HC_i(P/\C\{t\})$.
Note that there is an exact sequence
\[\xymatrix{
0 \ar[r] &\Omega^{2i-1}_f/d\Omega_f^{2i-2}\ar[r]^-{d}& 
\Omega^{2i}/(df\wedge d\Omega^{2i-1})
\ar[r]& \Omega^{2i}/d\Omega^{2i-1} \ar[r]& 0
}
\]

and $d\eta^{i-1}(M)=\omega^i(M)$ by definition.
If the number of variables $n+1$ is {\em even}, then a top degree form sits
in the {\em Brieskorn-module} 
\[{\mathcal H}^{(0)}_f:=\Omega^{n+1}/df \wedge d\Omega^{n-1}\,,\]
a free $\C\{t\}$-module of rank $=\mu(f)$, see \cite{looijenga}.\\
(Note however, that {\em only in the quasi-homogeneous case} its stalk at $0$,  
\[\mathcal{H}_f^{(0)}(0)=\Omega^{n+1}/(df \wedge d\Omega^{n-1}+f\Omega^{n+1})\,,\]
can be identified with the jacobian ring $P/J_f\cong \Omega^{n+1}/df \wedge \Omega^n$.)

The higher residue pairing
\[ K: \mathcal{H}_f^{(0)} \times\mathcal{H}_f^{(0)} \lra \C\{\!\{\partial_t^{-1}\}\!\}\]
of K.~Saito can be seen as the de Rham-realisation of the Seifert form of
the singularity. Expressing $lk(ch(M),ch(N))$ in this way
leads to a {\em Kapustin-Li type formula} for $\theta(M,N)$.
Details of this will appear elswhere.} 
\end{remark}

\section{ Behaviour under deformation}
As we have identified  the Theta pairing $\theta(M,N)$ as a topological
invariant, one would expect that this number remains {\em constant}, if
we let the hypersurface and modules depend on an additional parameter $t$.
We will now  indicate in what sense this is indeed the case.
As the arguments are well-known in singularity theory, we will be brief.
We consider a one-parameter deformation of an isolated hypersurface singularity
$f$:
\[F \in P:=\C\{x_0,x_1,\ldots,x_n,t\}\,,\quad F(x,0)=f(x)\,.\]
The hypersurface ring $R:=P/(F)$ is a flat module over the discrete
valuation ring  $\C\{t\}$, so we have an exact sequence
\[ 0 \lra R \stackrel{t \cdot}{\lra} R \lra \overline{R} \lra 0\,,\]
where $\overline{R}:=\C\{x_0,\ldots,x_n\}/(f)$. 
If $M$ and $N$ are $R$-modules that are flat as $\C\{t\}$-modules, they appear similarly in
exact sequences
\[ 0 \lra M \stackrel{t \cdot}{\lra} M \lra \overline{M} \lra 0\,,\quad
0 \lra N \stackrel{t \cdot}{\lra} N \lra \overline{N} \lra 0\,.
\]
From basic principles of homological algebra it follows that there is
a long exact sequence that in large enough (homological) degrees looks like
\[
\xymatrix@R7.5pt{
\cdots\ar[r]&Tor^R_{ev}(M,N)\ar[r]^-{t\cdot}&Tor^R_{ev}(M,N) \ar[r]&
Tor^{\overline{R}}_{ev}(\overline{M},\overline{N}) \ar[r]&\\
\ar[r]&Tor^R_{odd}(M,N)\ar[r]^-{t\cdot}&Tor^R_{odd}(M,N) \ar[r]&
Tor^{\overline{R}}_{odd}(\overline{M},\overline{N}) \ar[r]&\cdots
}
\]
If $length(Tor^{\overline{R}}_{ev/odd}(\overline{M},\overline{N})) <\infty$,
then the $R$-modules $Tor^R_{ev/odd}(M,N)$ are finitely generated as $\C\{t\}$-modules. 
Any finitely generated $\C\{t\}$-module $H$ is of the form $T\oplus F$, where $T$ is
the torsion submodule and $F$ is $\C\{t\}$-free of finite rank, $Rank(H)=Rank(F)$. In particular,
\[length(Coker(t\cdot :H \lra H))-length(Ker(t \cdot: H \lra H))=Rank(H)\,.\]
From this we immediately obtain\\

\begin{corollary}
\[\theta(\overline{M},\overline{N})=Rank(Tor^R_{even}(M,N))-Rank(Tor^R_{odd}(M,N))\,.\]
$\hfill\Diamond$
\end{corollary}
Now choose an appropriate representative $p: Y \lra S$ for the hypersurface $F=0$, where $S$ is a 
small disc in $\C$. We obtain coherent sheaves $Tor^{\calo_Y}_k(M,N)$ on 
$Y$. Under the same assumptions these sheaves are, via $p_*$, coherent as $\calo_S$-modules. The 
above $Rank(Tor^R_{even/odd}(M,N))$ can then be identified with
a sum of  contributions in a general fibre over $t \in \Delta$. Hence we obtain\\

\begin{theorem}
\[ 
\theta(\overline{M},\overline{N})=\sum_{y \in Y_t} \theta_{y}(M_{y},N_{y})\,,\]

where $\theta_{y}(M_{y},N_{y})$ denotes Hochster's Theta pairing over the local ring of 
$(Y_t,y)$ of the localizations $M_y$ and $N_y$ at $y \in Y_t$, and only finitely many 
summands are nonzero.$\hfill\Diamond$

\end{theorem}
Let us say that $M$ is a {\em smoothing} of $\overline{M}$ if it is locally 
free outside the fibre $t=0$, and call $\overline{M}$ {\em smoothable}, if it
has a smoothing.
\medskip

From the above theorem we then deduce:

\begin{corollary} If $\overline{M}$ is smoothable, then $\theta(\overline{M},\overline{N})=0$ 
for each $\overline{R}$-module $\overline{N}$ that extends to a $R$--module $N$ that is flat
as $\C\{t\}$-module. In particular, $\theta(\overline{M},\overline{M})=0$.$\hfill\Diamond$
\end{corollary}

To conclude, we revisit the examples 1.1. With $M=\calo_{L}=\C\{x,y\}/(x)$ the
maximal Cohen-Macaulay module associated to a component of the one-dimensional $A_{1}$
singularity $f=xy$, it is clear that this cycle cannot be extended to the family of hyperbolas 
$F=xy-t$, and the calculation $\theta(M,M)=1$ bears this out algebraically.

On the other hand, the line $L$ given by $x = z = 0$ on the two-dimensional $A_{1}$ 
singularity $f=xy-z^{2}$ lifts to a line in a ruling on the hyperboloid $F=xy-z^{2}-t$, explaining 
geometrically why the associated maximal Cohen-Macaulay module $M$ satisfies 
$\theta(M,M)=0$, as we saw.

\end{document}